\documentclass[preprint,12pt]{elsarticle}

\usepackage{amsmath,amssymb,amsthm}

\newtheorem{thm}{Theorem}
\newtheorem{lem}[thm]{Lemma}
\newdefinition{rmk}{Remark}
\newproof{pf}{Proof}
\newproof{pf1}{Proof of Theorem \ref{th:actsym}}
\newproof{pf2}{Proof of Theorem \ref{th:unique}}
\newproof{pf3}{Proof of Theorem \ref{th:actmin}}
\newproof{p-lem}{Proof of Lemma \ref{lem.weakMP}}

\journal{Journal of Differential Equations}

\begin{document}

\newcommand{\norm}[1]{\left\Vert#1\right\Vert}
\newcommand{\al}{\alpha}
\newcommand{\be}{\beta}
\newcommand{\de}{\partial}
\newcommand{\la}{\lambda}
\newcommand{\La}{\Lambda}
\newcommand{\ga}{\gamma}
\newcommand{\ep}{\epsilon}
\newcommand{\del}{\delta}
\newcommand{\Del}{\Delta}
\newcommand{\sig}{\sigma}
\newcommand{\ome}{\omega}
\newcommand{\Ome}{\Omega}
\newcommand{\uone}{u^{(1)}}
\newcommand{\utwo}{u^{(2)}}
\newcommand{\C}{{\mathbb C}}
\newcommand{\N}{{\mathbb N}}
\newcommand{\Z}{{\mathbb Z}}
\newcommand{\R}{{\mathbb R}}
\newcommand{\Rn}{{\mathbb R}^{n}}
\newcommand{\Rnu}{{\mathbb R}^{n+1}_{+}}
\newcommand{\Cn}{{\mathbb C}^{n}}
\newcommand{\spt}{\,\mathrm{supp}\,}
\newcommand{\Lin}{\mathcal{L}}
\newcommand{\SSS}{\mathcal{S}}
\newcommand{\F}{\mathcal{F}}
\newcommand{\xxi}{\langle\xi\rangle}
\newcommand{\xx}{\langle x\rangle}
\newcommand{\yy}{\langle y\rangle}
\newcommand{\dint}{\int\!\!\int}
\newcommand{\triple}[1]{{|\!|\!|#1|\!|\!|}}

\renewcommand{\Re}{\;\mathrm{Re}\;}

\begin{frontmatter}



\title{Symmetry  and uniqueness of minimizers of  Hartree type equations with external Coulomb potential}


\author{Vladimir Georgiev}

\address{Dipartimento di Matematica,  Universit\`a di Pisa\\
Largo Bruno Pintecorvo 5, 56127 Pisa,  Italy\\
e-mail: georgiev@dm.unipi.it}

\fntext[label2]{The first  author was supported by the Italian National Council of Scientific Research (project PRIN No.
2008BLM8BB )
entitled: "Analisi nello spazio delle fasi per E.D.P."}

\author{George Venkov}

\address{Faculty of Applied Mathematics and Informatics,
Technical University of Sofia\\
Kliment Ohridski 8,
1756 Sofia, Bulgaria\\
e-mail: gvenkov@tu-sofia.bg}

\begin{abstract}
In the present article we study the radial symmetry of minimizers
of the energy functional, corresponding to the repulsive Hartree
equation in external Coulomb potential. To overcome the
difficulties, resulting from the "bad" sign of the nonlocal term,
we modify the reflection method and then, by using Pohozaev
integral identities we get the symmetry result.
\end{abstract}


\end{frontmatter}


\section{Introduction}

Solitary waves associated with the Hartree type equation in
external Coulomb potential are solutions of type
$$
 \chi(x) e^{-i \omega t},\quad x\in \mathbb{R}^3 ,
t\in \mathbb{R},$$ where $\omega >0$ and $\chi $  satisfies the
nonlinear elliptic equation
\begin{eqnarray}\label{eq.M-S1}
 \ \ \ \ \ -\Delta \chi(x) + \int_{\mathbb{R}^3}\frac{|\chi(y)|^2 dy}{|x-y|} \ \chi(x) - \frac{\chi(x)}{|x|} + \omega \chi(x) =  0.
\end{eqnarray}

The natural energy functional associated with this problem is (see
\cite{CL1})
\begin{eqnarray}\label{eq.defEps}
    \mathcal{E} (\chi) = \frac{1}{2} \| \nabla \chi\|^2_{L^2}+ \frac{1}{4} A(|\chi|^2)- \frac{1}{2} \int_{\mathbb{R}^3}
    \frac{|\chi(x)|^2}{|x|}\ dx,
\end{eqnarray}
where we shall denote
\begin{equation}\label{eq.dAa}
   A(f) =  \int_{\mathbb{R}^3} \int_{\mathbb{R}^3}\frac{f(x)f(y)}{|x-y|} dy
    dx.
\end{equation}

The corresponding minimization problem is associated with the
quantity
\begin{equation}\label{defImu}
    I_N = \min \{ \mathcal{E}(\chi) ; \chi \in H^1, \|\chi\|^2_{L^2} = N
    \}.
\end{equation}

The existence of positive minimizers $\chi_0(x)$, such that
$$ \mathcal{E}(\chi_0) = I_N, \ \ \|\chi_0\|^2_{L^2} = N,$$
is established by Cazenave and Lions in \cite{CL1} by the aid of
the concentration compactness method.

For a given $\omega>0$, the constrained minimization problem
\eqref{defImu} can be compared with the unconstrained minimization
problem $$ S^{min}_\omega=\min\{S_\omega (\chi); \chi \in H^1 \},
$$where $S_\omega (\chi)$ is the corresponding action functional, defined by
\begin{equation}\label{eq.defSo}
   S_\omega (\chi) = \mathcal{E}(\chi) + \frac{\omega}{2}  \| \chi
   \|^2_{L^2}.
\end{equation}

There are  different results on the symmetry (and uniqueness) of the minimizers.
The basic result due to Gidas, Ni and Nirenberg \cite{GNN} implies
the radial symmetry of the minimizers associated with the
semilinear elliptic equation
$$ \Delta u + f(u) = 0,$$ provided suitable assumptions on the function $f(u)$
are satisfied and the scalar function $u$ is positive. As in the
previous result due to Serrin \cite{S}, the proof is based on the
maximum principle and the Hopf's lemma.

Therefore, the first natural question is to ask if the linear operator
$$ P_\omega = -\Delta  - \frac{1}{|x|} + \omega $$
in \eqref{eq.M-S1}, satisfies the weak maximum principle in the
sense that
\begin{equation}\label{eq.MPr}
    u \in H^2, \  P_\omega (u) =g \geq 0 ,  \Longrightarrow u \geq 0.
\end{equation}
The above maximum principle is incomplete, since additional
behavior of $ u$  and $g$ at infinity has to be imposed, namely,
we shall suppose that
\begin{equation}\label{mpdecay}
(1+|x|)^{-M}e^{\sqrt{\omega}|x|}u \in
   H^2,  \quad (1+|x|)^{-M}e^{\sqrt{\omega}|x|}g \in H^2,
\end{equation}
for some real number $M > 0.$

Note, that the energy levels of the hydrogen atom are described by
the eigenvalues $\omega_k > 0$ of the eigenvalue problem
$$ \Delta e_k(x) +   \frac{e_k(x)}{|x|} = \omega_k e_k(x), \ \ e_k(x) \in H^2. $$
One has
$$ \omega_k = \frac{1}{4(k+1)^2}, \quad k=0,1,... $$
and $e_0(x) = c e^{-|x|/2}, c >0.$ The first observation is that
all eigenfunctions $e_k(x)$, $k \geq 1$, are expressed in terms of
Laguerre polynomials of $|x|$, having exactly $k$ roots. This fact
guarantees that the maximum principle is not valid for $ \omega =
\omega_k.$ More precisely, we can show the following.

\begin{lem} \label{lem.weakMP} The weak maximum principle \eqref{eq.MPr} is valid if an only if
$$ \omega \geq \frac{1}{4}. $$
\end{lem}

This result can be compared with the existence of action
minimizers for the corresponding functional $S_\omega$, obtained
by Lions for $ 0 < \omega < 1/4$ (see for details \cite{PL}).

\begin{thm}\label{th:actmin} We have the properties:

a) for any $\omega > 0,$ the inequality $$ \min_{\chi \in H^1}
S_\omega (\chi) = S^{min}_\omega > - \infty$$ holds;

b) if $0 < \omega < 1/4$, then $ S^{min}_\omega <0;$

c) if $0 < \omega < 1/4$, then there exists a positive function
$\chi(x) \in H^1,$ such that
$$ S_\omega(\chi ) = S^{min}_\omega.$$
\end{thm}

Our main goal of this paper is to clarify if the positive
minimizers of $S_\omega$ are radially symmetric and unique. The
above results show that we have to consider the domain $ 0 <
\omega < 1/4,$ where the key tool of Gidas, Ni and Nirenberg (i.e.
the maximum principle for the corresponding linear operator) is
not applicable.

The symmetry of the energy functional (even with constraint
conditions) can not imply, in general, the radial symmetry of the
minimizers. This phenomena was discovered and studied in the works
\cite{CZE}, \cite{E1} and \cite{E2} in the scalar case.

Some sufficient conditions that guarantee the symmetry of
minimizers have been studied by Lopes in \cite{L}, by means of the
reflection method that (for the case of plane $x_1=0$) uses the
functions
\begin{equation*}
    u_1(x) = \left\{
       \begin{array}{ll}
         u(\hat{x}), \ \hat{x} =(-x_1,x_2,\cdots, x_n), & \hbox{if $x_1>0$;} \\
         u(x), & \hbox{if $x_1<0$}
       \end{array}
     \right.
\end{equation*}
and
\begin{equation*}
    u_2(x) = \left\{
       \begin{array}{ll}
         u(\hat{x}), \ \hat{x} =(-x_1,x_2,\cdots, x_n), & \hbox{if $x_1<0$;} \\
         u(x), & \hbox{if $x_1>0$.}
       \end{array}
     \right.
\end{equation*}

If the functional to be minimized has the form
$$
E(u) = \frac{1}{2} \| \nabla u\|^2_{L^2}+  \int_{\mathbb{R}^n}
    F(u(x))\ dx, $$
then we have the relation
$$ E(u_1)+E(u_2) = 2 E(u) $$
and this enables one to obtain the symmetry of minimizer, when
$F(u)$ is a combination of functions of type $|u|^p, p \geq 2$.

The reflection method works effectively when $u(x)$ is a
vector-valued function and constraint conditions (as in the
problem \eqref{defImu}) are involved too.

Recently, the reflection method was generalized in \cite{LM} and
\cite{M} for very general situations and one example of possible
application is the functional of type
$$ E(u) = \frac{1}{2} \| \nabla u\|^2_{L^2}+  \int_{\mathbb{R}^n}
    F(u(x))\ dx  - A(|u|^2), $$
involving nonlocal term as in \eqref{eq.defEps}. This Choquard
type functional has the specific property
$$
E(u_1)+E(u_2) \leq 2 E(u),$$ exploiting the negative sign of the
nonlocal term $A(|u|^2)$.

An analogous result for the scalar case can be obtained by means
of the Schwarz symmetrization (or spherical decreasing
rearrangement \cite{LL}) $u^*(|x|)$ of the non-negative $u \in
H^1.$ Indeed, we have the equality
$$  \int_{\mathbb{R}^n}
    F(u(x))\ dx  =  \int_{\mathbb{R}^n}
    F(u^*(x))\ dx, $$
as well as the inequalities
$$ \| \nabla u\|^2_{L^2} \geq \| \nabla u^*\|^2_{L^2}, \ \ \  A(|u|^2) \leq A(|u^*|^2), $$
so, we get
$$  E(u^*) \leq E(u)$$
and one can use the property that $u$ is minimizer.

The functional in \eqref{eq.defEps} is a typical example, when
reflection method and Schwarz symmetrization meet essential
difficulty to be applied directly.

The main goal of this work is to find an approach to establish the
symmetry of the minimizer for functionals of Hartree type
\eqref{eq.defEps}, involving nonlocal terms with "bad" sign.

To state this main result, we shall try  first to connect the
minimizers of the constraint problem \eqref{defImu} (associated
with the energy functional $\mathcal{E}(\chi )$) with the
minimization of the action functional $S_\omega(\chi )$. Similar
relation for local type interactions is discussed in chapter IX of
\cite{Caz}. Then, we shall establish that the minimizer of Theorem
\ref{th:actmin} is a radially symmetric function.

\begin{thm} \label{th:actsym} The solution $\chi(x)$ from Theorem \ref{th:actmin} is a radially symmetric function for
$$ \frac{1}{16} < \omega < \frac{1}{4}.$$
\end{thm}

\begin{rmk} The result of Theorem III.1 in \cite{CL1} treats more general  case of potentials of type
$$ V(x) = - \sum_{j=1}^K  \frac{Z}{|x-x_j|}, $$
while in our case we have
$$ V(x) = - \frac{Z}{|x|}.$$
Therefore,  the energy functional $ \mathcal{E} (\chi)  $ is
rotationally invariant in our case.  From Theorem \ref{th:actmin}
and Theorem \ref{th:actsym} one can see that the solution $ \chi_0
(x)$ of \eqref{defImu} is radially symmetric and unique (up to a
multiplication with complex number $z$, with $|z|=1$).
\end{rmk}
As it was mentioned above the energy (and therefore the action) is
a functional involving the nonlocal term with "bad" sign. To
explain the main idea to treat this case, we recall the rotational
symmetry of the energy (and action) functional. Therefore, if
$\chi$ is the action minimizer from Theorem \ref{th:actmin}, it is
sufficient to show that the solution is symmetric with respect to
$x_1$-plane, for any choice of the $x_1$-direction. In other
words, we consider $\hat{\chi}(x)=\chi(\hat{x})$, with
$\hat{x}=(-x_1,x_2,x_3)$ and we aim to prove that
$\chi=\hat{\chi}.$

To show this, we shall consider the two terms
$$
\chi_{\pm}=\frac{\chi \pm \hat{\chi}}{2}.$$

So, our goal is to verify the inequality
\begin{equation}\label{eq.mCl23}
   S_{\omega}(\chi_{+})+S_{\omega}(\chi_{-}) \leq
S_{\omega}(\chi)
\end{equation}
and see that the condition $\chi \neq \hat{\chi}$ implies
$S_{\omega}(\chi_{-}) >0.$

The form of the functional $S_{\omega}$ suggests one, in order to
verify \eqref{eq.mCl23}, to use  an appropriate version of the
Clarkson inequality for the quadratic form $A(f)$. Namely, we can
prove that the following inequality
$$ A\left(\left(\frac{f+g}{2}\right)^2\right)+A\left(\left(\frac{f-g}{2}\right)^2\right)\leq\frac{A(f^2)+A(g^2)}{2}$$
 holds true. Unfortunately, the usual Clarkson inequality in the form given above, is too rough to serve as a tool for proving \eqref{eq.mCl23}.
Therefore, we shall use a refined version of Clarkson inequality
(see Lemma \ref{Clarkson} below) in the form
\begin{eqnarray}
 \nonumber  A\left(\left(\frac{f+g}{2}\right)^2\right)+A\left(\left(\frac{f-g}{2}\right)^2\right)\leq\frac{A(f^2)+A(g^2)}{8}\\
 +\frac{3\sqrt{A(f^2)A(g^2)}}{4}.  \nonumber
\end{eqnarray}

The final step is to treat the uniqueness of positive minimizers.
 of
the problem
\begin{equation}\label{uniq.act-mini}
    S_{\omega}^{min} = \min \{ S_{\omega}(\chi) ; \chi \in H^1\}.
\end{equation}
 Our proof can not follow the Lieb's uniqueness proof
for the ground state solution of the Choquard equation \cite{Lie}.
In general, the Lieb's proof strongly depends on the specific
features of the nonlocal nonlinear equation \eqref{eq.M-S1} and
differs from the corresponding results for semilinear elliptic
equation given by Kwong in \cite{Kw}. Indeed, once the radial
symmetry is established, one can use Pohozaev identities and
reduce the nonlocal nonlinear elliptic problem \eqref{eq.M-S1} to
an ordinary differential equation of the type
$$ u''(r)+W(r) u(r) +4\pi\int^{r}_{0} \left(\frac{1}{s}-\frac{1}{r} \right)u^2(s)ds u(r)=\omega u(r),$$
where
$$ W_\chi (r) = \frac{1}{r} - 4\pi\int^{\infty}_{0} \chi^2(s)sds.$$
The positive sign in front of the nonlinear term is the main
obstacle to apply Sturm type argument and derive the uniqueness of
positive solutions to this ordinary differential equation.
However, for $\frac{1}{16} < \omega < \frac{1}{4}$ we can apply
the approach based on the refined Clarkson inequality and using
the orthogonal projection on the eigenspace of the first
eigenvalue of the operator $\Delta+1/|x|$, we can establish the
following result.
\begin{thm} \label{th:unique} Let $\frac{1}{16} < \omega < \frac{1}{4}$. Then, the solution $\chi$ of minimization problems \eqref{uniq.act-mini}  is unique.
\end{thm}

Let's mention that the results in Theorems \ref{th:actsym} and \ref{th:unique} can be compared with the results in
\cite{BS}, where the uniqueness of minimizers for the constrained variational problem \eqref{defImu} is studied.
To show the relations between action minimization and \eqref{defImu} one has to apply the uniqueness of action minimizers
or alternatively the uniqueness of minimizers of constrained variational problem.

The plan of the work is the following.  In Section 2 we consider
the maximum principle for the linear Schr\"odinger equation with
Coulomb potential and prove Lemma \ref{lem.weakMP}. The proof of
Theorem \ref{th:actsym}, stating that the minimizers are radially
symmetric is presented in Section 3 by the aid of a refined
version of Clarkson inequality. In Section 4 we establish the
Pohozaev integral relations, corresponding to equation
\eqref{eq.M-S1}, and in Section 5 we prove uniqueness Theorem
\ref{th:unique}. Finally, in \ref{app.A} we prove for completeness
the existence of positive action minimizers, stated in Theorem
\ref{th:actmin}, while in  \ref{app.C} the connection between
energy and action minimizers is discussed.

The authors are grateful to Louis Jeanjean  for important
discussions and remarks  on symmetry of minimizers as well as to the referee for pointing out a gap in the proof of the Theorem 3.

\section{Maximum principle for Schr\"odinger equation with Coulomb potential}

The maximum principle, stated in \eqref{eq.MPr} will be verified
by the aid of the substitution
$$ u = \varphi w, \varphi(x)=\varphi(|x|),$$
where $\varphi$ is a radial function, satisfying the property
\begin{equation}\label{eq.posf}
   -\Delta \varphi - \frac{\varphi}{|x|}+\omega \varphi = h(|x|) \geq 0.
\end{equation}
Our goal is to construct $\varphi$, so that $\varphi(|x|) > 0 .$
We have several possibilities, depending on $\omega.$ If $\omega
> 1/4$, we shall show that such a function exists and it is of type
\begin{equation}\label{eq.typ}
    \varphi(r) = e^{-\beta r} Q(r), \ \ \beta = \sqrt{\omega}, \ Q(r) = Ar^2+Br+C.
\end{equation}
If $\omega = 1/4$, then we can take simply $ \varphi(r) = e^{-
r/2}.$ If $0 < \omega < 1/4$, we shall see that a function
$\varphi$ of type \eqref{eq.typ} exists, but $\varphi(r)$ changes
the sign for $r>0.$ Hence, this function gives a counterexample,
showing that the weak maximum principle \eqref{eq.MPr} is not
fulfilled in this case.

Therefore, to complete the proof of Lemma \ref{lem.weakMP}, we
have to explain how the existence of positive $\varphi(r)$,
satisfying \eqref{eq.posf} will imply the weak maximum principle
and then to construct in different cases the function $Q(r)$ in
\eqref{eq.typ}, so that \eqref{eq.posf} is satisfied.

\begin{p-lem} After the substitution $u = \varphi w$, we have
$$ P_\omega(u) = -\varphi \Delta w - 2 \nabla \varphi \nabla w + P_\omega(\varphi) w =
\varphi \Delta w + 2 \nabla \varphi \nabla w + h w.$$ If
$\varphi(|x|)>0$, then we can write
$$ -\Delta w - \frac{2}{\varphi} \nabla \varphi \nabla w + \frac{h}{\varphi} w = \frac{g}{\varphi}.$$
Choosing $M=1,$ we see that $$ \frac{g}{\varphi} \in H^2,$$ so we
can apply the classical maximum principle (since $h \geq 0$) and
obtain $w \geq 0.$ This argument shows that the maximum principle
is fulfilled if the function $\varphi(r)$ satisfies inequality
\eqref{eq.posf} and its polynomial term $Q(r)>0$ for $r\geq 0$.

To construct $Q$, we substitute $\varphi(r) = e^{-\beta r} Q(r)$
into \eqref{eq.posf}  and find that
$$ e^{\beta r}  rh(r) = -(2B+C(-2\beta+1)) - (6 A+B(-4\beta+1))r + (6\beta-1)Ar^2 .$$
We take for simplicity $A=1$ and $$ B = C(\beta-1/2),\, C =
\frac{12}{(2\beta-1)(4\beta-1)}.$$ Then the condition $\beta >
1/2$ implies that
$$ B = \frac{6}{(4\beta-1)}>0,\, C = \frac{12}{(2\beta-1)(4\beta-1)} >0 $$ so $Q(r)>0$ and
$$e^{\beta r}  rh(r) = (6\beta-1)r^2 \geq 0.$$
This argument completes the proof of the weak maximum principle for $\omega > 1/4.$

If $1/16 < \omega < 1/4$, then we can take the same $A,B,C$ and see that
$$e^{\beta r}  rh(r) = (6\beta-1)r^2 \geq 0.$$ Since $A=1$ and $C < 0$ in this case, the function $Q(r)$ changes the sign.

Finally, if $0 < \omega < 1/16$, then we choose
$$ A=0, B = -1, C = \frac{1}{1/2-\beta}$$
and then
$$ Q(r) = \frac{2}{1-2\beta} - r, \ \ e^{\beta r}  rh(r) = (1-4\beta)r \geq 0.$$
Again, it is clear that $Q(r)$ changes the sign, and the proof of
the Lemma is completed.
\end{p-lem}

\section{Radial symmetry of action minimizers}

Even in the non-local case, the problem that action and energy
minimizers are nonnegative functions, is easy to be proved.
Indeed, if $\chi(x) \in H^1$ is a real-valued minimizer of the
functional
\begin{equation}\label{eq.action-form}
   S_{\omega}(\chi) = \frac{1}{2} \| \nabla \chi\|^2_{L^2}+ \frac{1}{4} A(\chi^2)- \frac{1}{2} \int_{\mathbb{R}^3}
    \frac{\chi(x)^2}{|x|}\ dx +
    \frac{\omega}{2}\|\chi\|^2_{L^2},
\end{equation}
then $|\chi(x)|$ satisfies the inequality
$$ \| \nabla |\chi|\|^2_{L^2} \leq \| \nabla \chi\|^2_{L^2},$$ as well as the identities
$$ A(|\chi|^2) = A(\chi^2), \ \ \int_{\mathbb{R}^3} \frac{|\chi(x)|^2}{|x|}\ dx  = \int_{\mathbb{R}^3} \frac{\chi(x)^2}{|x|}\ dx ,$$
so $|\chi(x)| \geq 0$ is also a minimizer of $S_\omega.$

Let us define the bilinear form
\begin{equation}\label{eq.Lomega}
   L_{\omega}(\chi,\psi)=\langle(-\Delta-\frac{1}{|x|}+\omega)\chi,\psi\rangle_{L^2},
   \quad \omega>0
\end{equation}
and the corresponding quadratic form
\begin{equation}\label{eq.L-quadr}
   L_{\omega}(\chi)=\langle(-\Delta-\frac{1}{|x|}+\omega)\chi,\chi\rangle_{L^2}.
\end{equation}

The quadratic form $A(\chi)$ defined in \eqref{eq.dAa} generates
the corresponding bilinear form
\begin{equation}\label{eq.A-form}
   A(\chi,\psi)=\int_{\mathbb{R}^3} \int_{\mathbb{R}^3}\frac{\chi(x)\psi(y)}{|x-y|} dy
    dx.
\end{equation}

Then, the action functional $S_{\omega}$ can be written as
\begin{equation}\label{eq.S-LA}
   S_{\omega}(\chi)=\frac{1}{2}L_{\omega}(\chi)+\frac{1}{4}
   A(\chi^2).
\end{equation}

Also, for any function $\chi$ we shall denote
$\hat{\chi}(x)=\chi(\hat{x})$, where $\hat{x}=(-x_1,x_2,x_3)$ for
any choice of our $x_1$-axis. It is easy to check that
\begin{equation}\label{eq.SandL}
   S_{\omega}(\chi)=S_{\omega}(\hat{\chi}), \quad
  L_{\omega}(\chi)=L_{\omega}(\hat{\chi}).
\end{equation}

With our next result, we shall establish Clarkson type
inequalities for the forms $A$ and $L_{\omega}$. In fact, we shall
prove the
 Lemma.

\begin{lem}\label{Clarkson} The following inequalities hold
\begin{equation}\label{eq.Clark-L}
   L_{\omega}\left(\frac{f+g}{2}\right)+L_{\omega}\left(\frac{f-g}{2}\right)=\frac{L_{\omega}(f)+L_{\omega}(g)}{2},
\end{equation}
\begin{eqnarray}\label{eq.Clark-A}
 \nonumber  A\left(\left(\frac{f+g}{2}\right)^2\right)+A\left(\left(\frac{f-g}{2}\right)^2\right)\leq\frac{A(f^2)+A(g^2)}{8}\\
 +\frac{3\sqrt{A(f^2)A(g^2)}}{4}.
\end{eqnarray}
\end{lem}

\begin{pf}
It is easy to verify the relation
\begin{eqnarray}\label{eq.Afg}
\nonumber
A\left(\left(\frac{f+g}{2}\right)^2\right)+A\left(\left(\frac{f-g}{2}\right)^2\right)\\
=\frac{1}{16}A(f^2+g^2+2fg)+\frac{1}{16}A(f^2+g^2-2fg).
\end{eqnarray}

Note that from $$ A(a+b)+A(a-b)=2A(a)+2A(b),$$ equality
\eqref{eq.Afg} becomes
\begin{eqnarray}\label{eq.Afg-ineq}
\nonumber
A\left(\left(\frac{f+g}{2}\right)^2\right)+A\left(\left(\frac{f-g}{2}\right)^2\right)
=\frac{1}{8}\left[A(f^2+g^2)+4A((fg)^2)\right]\\
\nonumber=\frac{1}{8}\left[A(f^2)+A(g^2)+2A(f^2,g^2)+4A((fg)^2)\right]\\
\leq \frac{A(f^2)+A(g^2)}{8}+\frac{3\sqrt{A(f^2)A(g^2)}}{4},
\end{eqnarray}
which proves \eqref{eq.Clark-A}. The first relation
\eqref{eq.Clark-L} in the Lemma, follows directly.
\end{pf}

The next result will play the crucial role in the present study.
We shall prove the following Lemma.

\begin{lem}\label{ClarksonII} If $L_{\omega}(f) = L_{\omega}(g)$ and $\mu, \nu \geq 0$ satisfy $2(\mu^2+\nu^2)=1,$ then
\begin{equation}\label{eq.Clark-Ln}
   L_{\omega}\left(\mu f + \nu g \right)+L_{\omega}\left(\mu f- \nu g\right)= L_{\omega}(f).
\end{equation}
If $A(f^2) = A(g^2)$ and $\mu, \nu \geq 0$ satisfy $2(\mu^2+\nu^2)=1,$ then
we have
\begin{equation}\label{eq.Clark-1LN}
   A\left(\left(\mu f+ \nu g\right)^2\right)+A\left(\left(\mu f- \nu g \right)^2\right)\leq A(f^2).
\end{equation}
\end{lem}

\begin{pf} Setting $\mu_1= 2\mu,$ $ \nu_1=2\nu$, we apply \eqref{eq.Clark-L} with $f,g$ replaced by $\mu_1 f$ and $\nu_1g$ respectively. Thus,
 we get
\begin{equation}\label{eq.Clark-Ln1}
L_{\omega}\left(\frac{\mu_1 f + \nu_1 g}{2} \right)+L_{\omega}\left(\frac{\mu_1 f - \nu_1 g}{2} \right)
= \frac{\mu_1^2L_{\omega}(f)+\nu_1^2L_\omega(g)}{2}.
\end{equation}
From $ L_{\omega}(f)= L_{\omega}(g)$ and $\mu_1^2+\nu_1^2=2$, we
complete the proof of \eqref{eq.Clark-Ln}.

Similarly, applying \eqref{eq.Clark-A} and the assumption $A(f^2)
= A(g^2)$, we find
\begin{eqnarray}\label{eq.Clark-mn1}
\nonumber  A\left(\left(\frac{\mu_1
f+\nu_1g}{2}\right)^2\right)+A\left(\left(\frac{\mu_1f-\nu_1g}{2}\right)^2\right)
\leq\frac{\mu_1^4+\nu_1^4 +6\mu_1^2\nu_1^2}{8}A(f^2)
\end{eqnarray}
or, equivalently
\begin{eqnarray}\label{eq.Clark-mn}
A\left(\left(\frac{\mu f+\nu
g}{2}\right)^2\right)+A\left(\left(\frac{\mu f-\nu
g}{2}\right)^2\right) \leq 2(\mu^4+\nu^4 +6\mu^2\nu^2)A(f^2)
\end{eqnarray}

Consider now the homogeneous quartic
polynomial
\begin{equation}\label{4-tic}
2(\mu^4+\nu^4
+6\mu^2\nu^2)
\end{equation}
on the circle $\mu^2+\nu^2=\frac{1}{2}$. Substituting
$\nu^2=\frac{1}{2}-\mu^2$, we obtain the following estimate
\begin{eqnarray}\label{eq.polynom}
\nonumber 2(\mu^4+\nu^4 +6\mu^2\nu^2)=2((\mu^2+\nu^2)^2
+4\mu^2\nu^2)\\=\frac{1}{2} +4\mu^2 -8\mu^4
=1-\frac{(1-4\mu^2)^2}{2}\leq 1.
\end{eqnarray}

Then, from \eqref{eq.Clark-mn} and \eqref{eq.polynom} follows the
proof of the Lemma.
\end{pf}

Turning back to the minimization problem of the action functional
$S_\omega$, we observe the following fact. If $\chi(x)$ is a
minimizer of the problem
\begin{equation}\label{act-inf}
 \min_{\chi \in H^1}
S_\omega (\chi),
\end{equation}
then $\hat{\chi}(x)$  and $-\hat{\chi}(x)$ are also minimizers of
$S_{\omega}(\chi) $. Moreover, we have the property.

\begin{lem}\label{L-argument} Assume that $\chi(x)$ is a minimizer of the problem \eqref{act-inf} and one of the
following alternatives:
\begin{enumerate}
    \item $L_{\omega}(\chi - \hat{\chi})\geq 0$;
    \item $L_{\omega}(\chi + \hat{\chi})\geq 0$
\end{enumerate}
holds. Then $\chi = \hat{\chi}$.
\end{lem}

\begin{pf} For simplicity, we shall consider the first case
only. Suppose $\chi \neq\hat{\chi}$, then from \eqref{eq.Clark-L}
we have
\begin{equation}\label{eq.L-sum}
   L_{\omega}\left(\frac{\chi + \hat{\chi}}{2}\right)+L_{\omega}\left(\frac{\chi -
   \hat{\chi}}{2}\right)=L_{\omega}(\chi),
\end{equation}
implying
\begin{equation}\label{eq.L-sum1}
   L_{\omega}\left(\frac{\chi + \hat{\chi}}{2}\right)\leq L_{\omega}(\chi).
\end{equation}

On the other hand, it is easy to check that the following Cauchy
inequalities
\begin{equation}\label{eq.CauchyIneq}
   A(f^2,g^2)\leq \sqrt{A(f^2)A(g^2)}, \quad A(fg)\leq \sqrt{A(f^2)A(g^2)}
\end{equation}
hold true. Applying now \eqref{eq.Clark-A}, we obtain
\begin{eqnarray}\label{eq.A-sum}
 \nonumber  A\left(\left(\frac{\chi + \hat{\chi}}{2}\right)^2\right)+A\left(\left(\frac{\chi - \hat{\chi}}{2}\right)^2\right)\leq
 \frac{A(\chi^2)+A(\hat{\chi}^2)}{8}\\
 +\frac{3\sqrt{A(\chi^2)A(\hat{\chi}^2)}}{4}\leq\frac{A(\chi^2)+A(\hat{\chi}^2)}{2}=A(\chi^2),
\end{eqnarray}
which, together with the assumption $\chi \neq\hat{\chi}$ gives
that
\begin{equation}\label{eq.A-sum1}
   A\left(\left(\frac{\chi + \hat{\chi}}{2}\right)^2\right)<
   A(\chi^2).
\end{equation}

Thus, from \eqref{eq.L-sum1}, \eqref{eq.A-sum1} and the definition
\eqref{eq.S-LA} it follows
\begin{equation}\label{eq.S-sum}
   S_{\omega}\left(\frac{\chi + \hat{\chi}}{2}\right)\leq S_{\omega}(\chi),
\end{equation}
which  contradicts to the assumption that $\chi$ is a minimizer.
This proves the Lemma.
\end{pf}

Now, we are ready to prove the radial symmetry of the action
minimizer, stated in Theorem \ref{th:actsym}.

\begin{pf1}

Taking into account Lemma \ref{L-argument}, we shall take a minimizer $\chi(x) \geq 0$ of $S_\omega$ and
shall  show that the condition
$$ \frac{1}{16} < \omega < \frac{1}{4}, $$
implies
 that $\chi = \hat{\chi}$ or
\begin{equation}\label{eq.L-negative}
    L_{\omega}\left(\chi -
   \hat{\chi}\right) > 0.
\end{equation}
Let
$$ \chi(x) = e_0(x) + f(x),$$
where $e_0(x) = ce^{-|x|/2}, c >0$ is the eigenvector corresponding to the first eigenvalue of the operator
$\Delta+1/|x|,$ while $ \langle f,e_0 \rangle_{L^2}=0.$ Since $e_0$ is a radial function, we have
$\hat{e_0} = e_0,$ so
$$ \chi - \hat{\chi} = f - \hat{f} =g, \ \ \langle g, e_0 \rangle_{L^2}=0, \ g \neq 0.$$

\begin{lem}\label{gort} Let us assume that $g \perp e_0$ in $L^2$. Then
$$ L_{\omega}\left(g\right) \geq \left( \omega - \frac{1}{16}\right) \|g\|^2_{L^2}. $$
\end{lem}

\begin{pf}
Note that $g \perp e_0$ in $L^2$ implies
$$ g = \sum_{k \geq 1} c_k e_k + h, $$
where $h$ is in the absolutely continuous space of the self-adjoint operator $ \Delta+ \frac{1}{|x|}$ in $L^2$, while  $e_k$ are eigenvectors of the same operator in $\{g \in L^2;  g \perp e_0 \}$ with eigenvalues
$\omega_k \leq 1/16.$ On the absolutely continuous space the operator has spectrum on $(-\infty,0)$ and it is non positive, so
$$ \left\langle \left(\Delta+ \frac{1}{|x|} \right)  h, h \right\rangle  \leq 0. $$

Hence, we have
$$ \left\langle \left(\Delta+ \frac{1}{|x|} \right)  g, g \right\rangle \leq \sum |c_k|^2 \omega_k \leq \frac{1}{16}\left( \sum |c_k|^2 \right) = \frac{1}{16} \|g\|^2_{L^2}$$
and
$$  L_{\omega}\left(g\right) = - \left\langle \left(\Delta+ \frac{1}{|x|} \right)  g, g \right\rangle + \omega \|g\|^2_{L^2} \geq \left( \omega - \frac{1}{16}\right) \|g\|^2_{L^2}.$$
This completes the proof of the Lemma.
   \end{pf}
Applying the above Lemma, we find
$$  L_{\omega}\left(\chi -
   \hat{\chi}\right) = L_{\omega}\left(g\right) \geq  \left( \omega - \frac{1}{16}\right) \|g\|^2_{L^2} > 0,$$
since $\omega > 1/16$ and $g \neq 0.$
Hence, \eqref{eq.L-negative} is fulfilled and the proof of the Theorem is complete.

\end{pf1}

\section{Pohozaev identities}
\label{app.B}

In this part we shall establish the so-called Pohozaev identities
for \eqref{eq.M-S1}. More precisely, we shall prove the following

\begin{lem}\label{Pohoz} If $\chi\in H^1(\mathbb{R}^3)$ and
satisfies \eqref{eq.M-S1} in $ H^{-1}(\mathbb{R}^3)$, then the
following identities hold
\begin{equation}\label{eq.Poh-1}
   \| \nabla \chi\|^2_{L^2} +\omega \| \chi\|^2_{L^2}=  \int_{\mathbb{R}^3}
    \frac{|\chi(x)|^2}{|x|}\ dx-A(|\chi|^2),
\end{equation}

\begin{equation}\label{eq.Poh-2}
    \| \nabla \chi\|^2_{L^2} +3\omega \| \chi\|^2_{L^2}= 2 \int_{\mathbb{R}^3}
    \frac{|\chi(x)|^2}{|x|}\ dx-\frac{5}{2}A(|\chi|^2).
\end{equation}
\end{lem}

\begin{pf}
To prove \eqref{eq.Poh-1} we multiply equation \eqref{eq.M-S1} by
$\bar{\chi}$, take the real part and integrate over
$\mathbb{R}^3$. To prove \eqref{eq.Poh-2} we shall use the
following relations
\begin{equation}\label{eq.div-1}
    \nabla\cdot(x|\chi|^2)=3|\chi|^2+2\Re\chi(x\cdot
    \nabla\bar{\chi}),
\end{equation}
\begin{equation}\label{eq.div-2}
    \nabla\cdot\left(x|\nabla\chi|^2-2\Re\nabla\chi(x\cdot
    \nabla\bar{\chi})\right)=|\nabla\chi|^2-2\Re\Delta\chi(x\cdot
    \nabla\bar{\chi}),
\end{equation}
\begin{equation}\label{eq.div-3}
    \nabla\cdot(x\frac{|\chi|^2}{|x|})=2\frac{|\chi|^2}{|x|}+2\Re\frac{\chi(x\cdot
    \nabla\bar{\chi})}{|x|},
\end{equation}
and
\begin{eqnarray}\label{eq.div-4}
 \nonumber \nabla\cdot\left(x\int_{\mathbb{R}^3}\frac{|\chi(y)|^2 dy}{|x-y|}|\chi|^2\right)=3\int_{\mathbb{R}^3}\frac{|\chi(y)|^2
   dy}{|x-y|}|\chi|^2\\
     \qquad-\int_{\mathbb{R}^3}\frac{x(x-y)|\chi(y)|^2 dy}{|x-y|^3}|\chi|^2+2\int_{\mathbb{R}^3}\frac{|\chi(y)|^2 dy}{|x-y|}\Re\chi(x\cdot
    \nabla\bar{\chi}).
\end{eqnarray}

Integrating \eqref{eq.div-1}--\eqref{eq.div-4} over $\mathbb{R}^3$
implies the equalities
\begin{equation}\label{eq.ident-1}
    \Re\int_{\mathbb{R}^3}\chi(x\cdot \nabla\bar{\chi})\ dx=-\frac{3}{2} \| \chi\|^2_{L^2},
\end{equation}
\begin{equation}\label{eq.ident-2}
    \Re\int_{\mathbb{R}^3}\Delta\chi(x\cdot \nabla\bar{\chi})\ dx=\frac{1}{2} \| \nabla \chi\|^2_{L^2},
\end{equation}
\begin{equation}\label{eq.ident-3}
    \Re\int_{\mathbb{R}^3}\frac{1}{|x|}\chi(x\cdot \nabla\bar{\chi})\ dx=-\int_{\mathbb{R}^3}
    \frac{|\chi(x)|^2}{|x|}\ dx,
\end{equation}
\begin{eqnarray}\label{eq.ident-4}
   \nonumber \Re\int_{\mathbb{R}^3}\int_{\mathbb{R}^3}\frac{|\chi(y)|^2 \chi(x)(x\cdot \nabla\bar{\chi}(x))}{|x-y|} dydx=
   -\frac{3}{2}A(|\chi|^2) \\+ \frac{1}{2}\int_{\mathbb{R}^3}\int_{\mathbb{R}^3}
    \frac{x(x-y)|\chi(y)|^2|\chi(x)|^2}{|x-y|^3} dydx.
\end{eqnarray}

On the other hand, observing the symmetry
\begin{eqnarray}\label{eq.sym-1}
 \nonumber\int_{\mathbb{R}^3}\int_{\mathbb{R}^3}
    \frac{x(x-y)|\chi(y)|^2|\chi(x)|^2}{|x-y|^3}\ dydx\\=\int_{\mathbb{R}^3}\int_{\mathbb{R}^3}
    \frac{y(y-x)|\chi(y)|^2|\chi(x)|^2}{|x-y|^3}\ dydx,
\end{eqnarray}
we calculate
\begin{eqnarray}\label{eq.sym-2}
 \nonumber\int_{\mathbb{R}^3}\int_{\mathbb{R}^3}
    \frac{x(x-y)|\chi(y)|^2|\chi(x)|^2}{|x-y|^3}\ dydx\\=\frac{1}{2}\int_{\mathbb{R}^3}\int_{\mathbb{R}^3}
    \frac{(x-y)^2|\chi(y)|^2|\chi(x)|^2}{|x-y|^3}\
    dydx=\frac{1}{2}A(|\chi|^2).
\end{eqnarray}

Substituting \eqref{eq.sym-2} into \eqref{eq.ident-4} we get
\begin{eqnarray}\label{eq.id-4-new}
   \Re\int_{\mathbb{R}^3}\int_{\mathbb{R}^3}\frac{|\chi(y)|^2 \chi(x)(x\cdot \nabla\bar{\chi}(x))}{|x-y|} dydx=
   -\frac{5}{4}A(|\chi|^2).
\end{eqnarray}

Finally, multiplying equation \eqref{eq.M-S1} by $x\cdot
\nabla\bar{\chi}$, taking the real part, integrating over
$\mathbb{R}^3$ and using \eqref{eq.ident-1}, \eqref{eq.ident-2},
\eqref{eq.ident-3} and \eqref{eq.id-4-new} we complete the proof
of the Lemma.
\end{pf}

The  Pohozaev identities are useful to treat the uniqueness of the
minimizers (modulo multiplication by complex constant $z$ with
$|z|=1$). Indeed, let $\chi_1$ and $\chi_2$ are minimizers of the
problem
\begin{equation}\label{act-min}
    S_{\omega}^{min} = \min \{ S_{\omega}(\chi) ; \chi \in H^1\}.
\end{equation}

Since
$$ S_\omega (\chi) = \frac{1}{2} \| \nabla \chi\|^2_{L^2} + \frac{\omega}{2} \|\chi\|^2_{L^2}
- \frac{1}{2} \int_{\mathbb{R}^3} \frac{|\chi(x)|^2}{|x|}\ dx +
\frac{1}{4} A(|\chi|^2),$$ we can apply the Pohozaev identities of
Lemma \ref{Pohoz}. In this way  we find
\begin{equation}\label{eq.Sequ}
 S_\omega (\chi) = -\frac{1}{4} A(|\chi|^2)
\end{equation}
and
\begin{equation}\label{eq.uniqP}
 A(|\chi_1|^2) =  A(|\chi_2|^2), \ L_{\omega}(\chi_1) = L_{\omega}(\chi_2),
\end{equation}
where $L_{\omega}(\chi)$ is defined according to
\eqref{eq.L-quadr}.

\section{Uniqueness of minimizers}\label{app.D}

In this section we shall prove the uniqueness result of Theorem
\ref{th:unique}. The classical approach for proving the uniqueness
of minimizers is to reduce the initial nonlinear equation to an
ordinary differential equation, using the radial symmetry.
Uniqueness of positive ground state solutions for nonlinear
Schr\"odinger equation on $\R^n$ with local nonlinearities of the
form $|u|^pu$ for $0 < p < \frac{4}{n-2}$, is a well-known fact,
due to Kwong \cite{Kw}. The proof in this case relies on Sturm
comparison theorems, but it cannot be applied directly to nonlocal
equations, such as \eqref{eq.M-S1}. For the attractive Choquard
equation, Lieb in \cite{Lie} prove uniqueness of energy minimizer
by using Newton's theorem for radial function $f(x)=f(|x|)$, that
is
\begin{eqnarray}\label{eq.maxgap-pr}
\int \frac{f(|y|)}{|x-y|^{n-2}} dy =   \int \frac{f(|y|)}{\max
  \{|x|,|y|\}}dy.
\end{eqnarray}.

The repulsive sign of the Hartree term in \eqref{eq.M-S1} is again
the main obstacle for applying directly the standard technique.

\begin{pf2}

Let $\chi_1$ and $\chi_2$ are non negative minimizers of the problem
$$
    S_{\omega}^{min} = \min \{ S_{\omega}(\chi) ; \chi \in H^1\}.
$$
Since they are radial functions, one can rewrite the elliptic
equation \eqref{eq.M-S1}, using Newton's theorem
\eqref{eq.maxgap-pr}, as an ordinary differential equation of the
form
\begin{equation}\label{eq.ODE}
-\chi''(r)-\frac{2}{r}\chi'(r)- \frac{\chi(r)}{r}
+4\pi\int^{\infty}_{0} \frac{\chi^2(s)s^2ds}{\max
  \{r,s\}}\chi(r)+ \omega \chi(r)=0.
\end{equation}

The above equation can be rewritten in the form
$$
-\chi''(r)-\frac{2}{r}\chi'(r)- W(r) \chi(r) +4\pi\int^{r}_{0}
\chi^2(s)\left(\frac{1}{r}-\frac{1}{s} \right)s^2ds \chi(r)+
\omega \chi(r)=0,
$$
where
$$ W(r) = \frac{1}{r} - 4\pi\int^{\infty}_{0} \chi^2(s)sds.$$

If we set $u(r) = r\chi(r)$, then from the identity
$$
\chi''(r)+\frac{2}{r}\chi'(r)= \frac{
 u''(r)}{r},
$$
the last equation becomes
\begin{equation}\label{eq.ODE-u}
u''(r)+W(r) u(r) -4\pi\int^{r}_{0} \left(\frac{1}{r}-\frac{1}{s} \right)u^2(s)ds u(r)=\omega u(r).
\end{equation}
This observation shows that the assumption $\chi(x)$ is a non negative minimizer implies $u(r)>0$ for $r>0.$
Hence $\chi_1(x)$ and $\chi_2(x)$ are positive functions.

Our goal is to use the projection of $\chi_1$ and $\chi_2$ on the one dimensional eigenspace
$$ E_0= \{ \alpha e^{-|x|/2}, \alpha \in (-\infty, \infty) \}$$ is the eigenvector corresponding to the first eigenvalue $\omega_0 = 1/4$ of the operator
$\Delta+1/|x|.$ First, we have to observe that  $\chi_1$  is not
orthogonal to $E_0.$ Indeed, if $\chi_1 \perp E_0$, then Lemma
\ref{gort} implies
$$L_\omega(\chi_1) \geq \left( \omega - \frac{1}{16} \right) \|\chi_1\|^2_{L^2} >0.$$
The relation \eqref{eq.S-LA} guarantees now $S_{\omega}(\chi_1)>0$ and this contradicts the relation \eqref{eq.Sequ}.
The contradiction shows that $ \chi_1 $ (and also $\chi_2$) is not orthogonal to $E_0.$

Let $$ \chi_1 = \mu_1 \alpha e^{-|x|/2} + f_1,\ \ \chi_2 = \mu_2
\alpha e^{-|x|/2} + f_2, $$ where $\alpha e^{-|x|/2} \in E_0,$
with $\alpha > 0$ and $ f_1, f_2 \perp E_0.$ Note that $\mu_1,
\mu_2 > 0$, since $\chi_1, \chi_2$ and $e_0$ are positive
functions. We  can choose  $\alpha>0$, such that
\begin{equation}\label{eq.relmn}
    2(\mu_1^2 + \mu_2^2)=1,
\end{equation}
used as assumption in Lemma \ref{ClarksonII}.
The other assumption
$$
 A(|\chi_1|^2) =  A(|\chi_2|^2), \ L_{\omega}(\chi_1) = L_{\omega}(\chi_2),
$$
is already established in \eqref{eq.uniqP}.

Applying Lemma \ref{ClarksonII}, we find the identity
$$ L_{\omega}\left(\mu_2 \chi_1 + \mu_1 \chi_2 \right)+L_{\omega}\left(\mu_2 \chi_1 - \mu_1 \chi_2\right)= L_{\omega}(\chi_1),$$
as well as the
inequality
$$ A\left(\left(\mu_2 \chi_1 + \mu_1 \chi_2\right)^2\right)+A\left(\left(\mu_2 \chi_1 - \mu_1 \chi_2\right)^2\right)\leq  A(\chi_1^2).$$

Then, we have the relation
$$ \mu_2 \chi_1 - \mu_1 \chi_2 = \mu_2 f_1 - \mu_1 f_2 = g \perp E_0.$$
If $g=0$, then $\chi_1= \mu_1 \chi_2/\mu_2$ and one can use the
ODE \eqref{eq.ODE-u} and the corresponding integral identities
\eqref{eq.Poh-1} and \eqref{eq.Poh-2}, to show that
$\chi_1=\chi_2.$ If $g \neq 0,$ then one can apply Lemma
\ref{gort} and find
$$ L_{\omega}\left(\mu_2 \chi_1 - \mu_1 \chi_2\right) \geq \left( \omega - \frac{1}{16} \right) \|g\|^2_{L^2} >0.$$
Hence,
$$ S( \mu_2 \chi_1 + \mu_1 \chi_2) = \frac{1}{2}L_{\omega}(\mu_2 \chi_1 + \mu_1 \chi_2)+\frac{1}{4}
   A((\mu_2 \chi_1 + \mu_1 \chi_2)^2) < S_\omega(\chi_1)$$
and this is a contradiction. The contradiction shows that
$\chi_1=\chi_2$ and this completes the proof of  Theorem
\ref{th:unique}.
\end{pf2}

\appendix

\section{Existence of action minimizers}\label{app.A}

The existence of action minimizers for Hartree type equation is
already established in \cite{PL}. For completeness, we shall
sketch the proof.

To show the boundedness from below of $S_\omega$, we shall prove
the following inequalities involving homogeneous Sobolev norms
$$ \|f\|_{\dot{H}^s(\mathbb{R}^3)} = \|(-\Delta)^{s/2} f\|_{L^2( \mathbb{R}^3)},  \ \ s > - 3/2.$$

\begin{lem} \label{lem.in1}For any $ p_1 \in [3,6]$ and $ p_2 \in [2,3]$ we have the estimates
\begin{equation}\label{eq.Ex1}
  \left( \int_{|x| \leq 1} |\chi(x)|^{p_1} dx\right)^{1/p_1} \leq C \| \chi \|^{\theta_1}_{\dot{H}^1}
  \| \chi^2 \|^{(1-\theta_1)/2}_{\dot{H}^{-1}}
\end{equation}
\begin{equation}\label{eq.Ex2}
  \left( \int_{|x| \geq 1} |\chi(x)|^{p_2} dx\right)^{1/p_2} \leq C \|\chi\|^{\theta_2}_{L^2}\| \chi \|^{\theta_3}_{\dot{H}^1}
  \| \chi^2 \|^{(1-\theta_2-\theta_3)/2}_{\dot{H}^{-1}},
\end{equation}
where
$$ \theta_1 = \frac{5}{3} - \frac{4}{p_1},\ \  \theta_2 = \frac{4(3-p_2)}{p_2}, \ \ \theta_3 = \frac{p_2-2}{p_2} . $$
\end{lem}

\begin{rmk} The assumptions $ p_1 \in [3,6]$ and $ p_2 \in [2,3]$ guarantee that all parameters
$ \theta_1, \theta_2, \theta_3, \theta_2+\theta_3 $ are in the
interval $[0,1].$
\end{rmk}
\begin{rmk} The relation
$$ \| f \|^2_{\dot{H}^{-1}} = \langle (-\Delta)^{-1} f, f \rangle_{L^2} = \frac{1}{4\pi}
\int_{\mathbb{R}^3} \int_{\mathbb{R}^3}\frac{f(x)f(y)}{|x-y|} dy
dx$$ implies
$$  \| \chi^2 \|^2_{\dot{H}^{-1}} = \frac{1}{4\pi} A(\chi^2).$$
\end{rmk}

\begin{pf} For $p_1=6$ the inequality \eqref{eq.Ex1} becomes
$$ \left( \int_{|x| \leq 1} |\chi(x)|^{6} dx\right)^{1/6} \leq C \| \chi \|_{\dot{H}^1}
  $$
and this is the standard Sobolev embedding. For $p_1=3$ we have to
verify the following estimate
\begin{equation}\label{eq.Ex3}
    \left( \int_{\mathbb{R}^3} |\chi(x)|^{3} dx\right)^{1/3} \leq C \| \chi \|^{1/3}_{\dot{H}^1}
  \| \chi^2 \|^{1/3}_{\dot{H}^{-1}}.
\end{equation}
  This inequality follows from
  $$\left| \int f(x) g(x) dx \right| \leq \| f \|_{\dot{H}^1} \| g \|_{\dot{H}^{-1}}$$
  with $f(x) = |\chi(x)|, $ $g(x) = |\chi(x)|^2 = \chi^2(x)$ and the observation that
 $$ \| |\chi| \|_{\dot{H}^1} = \| \chi \|_{\dot{H}^1}.$$
 Interpolation between $p_1=6$ and $p_1=3$ proves \eqref{eq.Ex1}.

 The inequality \eqref{eq.Ex2} for $p_2=3$ follows from  \eqref{eq.Ex3}.

 For $p_2=2$  \eqref{eq.Ex2} reduces to the simple inequality
 $$ \left( \int_{|x| \geq 1} |\chi(x)|^{2} dx\right)^{1/2} \leq C \|\chi\|_{L^2}. $$
 An interpolation argument implies \eqref{eq.Ex2} and completes the proof of the Lemma.

\end{pf}

After this Lemma we can show that the action functional is bounded
from below.

\begin{lem} \label{lem:bound} For any $\omega > 0 $ the inequality $$ \min_{\chi \in H^1} S_\omega (\chi) = S^{min}_\omega > - \infty$$ holds. For
 $0 < \omega < 1/4$ we have $ S^{min}_\omega <0.$

\end{lem}

\begin{pf} The only negative term in $S_\omega$ is
$$ - \frac{1}{2} \int_{\mathbb{R}^3}
    \frac{|\chi(x)|^2}{|x|}\ dx. $$ Decomposing the integration domain into $ |x| \leq 1 $ and $ |x| > 1$ we apply H\"older inequality and obtain
    $$  \int_{\mathbb{R}^3} \frac{|\chi(x)|^2}{|x|}\ dx \leq C \left( \int_{|x| \leq 1} |\chi(x)|^{p_1} dx\right)^{2/p_1}
    + C \left( \int_{|x| > 1} |\chi(x)|^{p_2} dx\right)^{2/p_2}, $$
    where $p_1 > 3 > p_2$. Applying Lemma \ref{lem.in1} as well as the Young inequality
    $$ X^{\theta_1} Y^{\theta_2} Z^{\theta_3} \leq \varepsilon X + \varepsilon Y + C_\varepsilon Z,$$
    with
    $$\theta_j \in (0,1), \theta_1+\theta_2+\theta_3 = 1,$$ we get
$$  \int_{\mathbb{R}^3} \frac{|\chi(x)|^2}{|x|}\ dx \leq \varepsilon \|\chi\|^{2}_{L^2} + \varepsilon \|\nabla \chi\|^{2}_{L^2}
+ C_\varepsilon \sqrt{A(\chi^2)}.
 $$
This estimate implies
$$ S_\omega(\chi) \geq \frac{1-\varepsilon}{2} \| \nabla \chi\|^2_{L^2}+ \frac{\omega-\varepsilon}{2} \| \chi\|^2_{L^2} + \frac{1}{4} A(\chi^2)- C_\varepsilon \sqrt{A(\chi^2)}. $$
Choosing $\varepsilon >0$ so small that $\varepsilon < \min (1,
\omega),$ we find
$$ S_\omega(\chi) \geq  \frac{1}{4} A(\chi^2)- C_\varepsilon \sqrt{A(\chi^2)} \geq - 2C_\varepsilon^2.$$

To finish the proof we take $\chi_\delta(x) = \delta e^{-|x|/2},$
such that
$$\left(\Delta+ \frac{1}{|x|} \right) \chi_\delta = \frac{1}{4} \chi_\delta.$$
Then
$$ 2S_\omega(\chi_\delta) = (\omega - 1/4) \|\chi_\delta\|^2_{L^2} + A(\chi_\delta^2)/2.$$
Since
$$ \|\chi_\delta\|^2_{L^2} = C_0 \delta^2, \ \ A(\chi_\delta^2)/2 = O(\delta^4),$$
the condition $\omega \in (0, 1/4)$ implies $
2S_\omega(\chi_\delta)  <0$ and this completes the proof.
\end{pf}

\begin{pf3} Take a minimizing sequence $\chi_k \in H^1,$ so that
\begin{equation}\label{eq.prs1}
    \lim_{k \rightarrow \infty }S_\omega(\chi_k) = S_\omega^{min} < 0.
\end{equation}
The argument of the proof of Lemma \ref{lem:bound} guarantees that
there exists a constant $C>0,$ so that
\begin{equation}\label{eq.bh1}
   \| \chi_k \|_{H^1} \leq C.
\end{equation}
One can find $\chi_*(x) \in H^1$ so that (after taking a
subsequence) $\chi_k$ tends weakly in $H^1$ to $\chi_*.$ Using the
inequality
$$
\int_{|x|>R} \frac{|\chi(x)|^2}{|x|} \ dx \leq \frac{C}{R},$$ and
the compactness of the embedding $L^p(|x|<R) \hookrightarrow
H^1(|x|<R),$ when $2 \leq p < 6,$ we see that (choosing a suitable
subsequence)
\begin{equation}\label{eq.bhh2}
    \lim_{k \rightarrow \infty} \int_{\mathbb{R}^3} \frac{|\chi_k(x)|^2}{|x|} \ dx =
\int_{\mathbb{R}^3} \frac{|\chi_*(x)|^2}{|x|} \ dx .
\end{equation}

Then we introduce $\varphi_k, \ \varphi_*$ so that
$$ \Delta \varphi_k = - 4\pi \chi_k^2(x),  \ \ \Delta \varphi_* = - 4\pi \chi_*^2(x).$$
One can show that $ \varphi_k$ tends weakly to $\varphi_*$ in
$\dot{H}^1.$ We have also the identities
$$ A(\chi_k^2) = \int \varphi_k(x) \chi^2_k(x) dx =  \frac{1}{4\pi} \| \nabla \varphi_k \|^2_{L^2} $$
and
$$ A(\chi_*^2) = \int \varphi_*(x) \chi^2_*(x) dx =  \frac{1}{4\pi} \| \nabla \varphi_* \|^2_{L^2} $$
so we obtain
$$ S_\omega(\chi_k) = \frac{1}{2} \| \nabla \chi_k\|^2_{L^2}+ \frac{\omega}{2} \| \chi_k\|^2_{L^2} + \frac{1}{4} A(\chi_k^2)-\frac{1}{2}
 \int_{\mathbb{R}^3} \frac{|\chi_k(x)|^2}{|x|} \ dx   $$
 $$ =\frac{1}{2} \| \nabla \chi_k\|^2_{L^2}+ \frac{\omega}{2} \| \chi_k\|^2_{L^2} + \frac{1}{16\pi} \| \nabla \varphi_k \|^2_{L^2}-
\frac{1}{2} \int_{\mathbb{R}^3} \frac{|\chi_k(x)|^2}{|x|} \ dx.$$

 Using \eqref{eq.prs1} and \eqref{eq.bhh2}, we get
$$ \lim_{k \rightarrow \infty} S_\omega(\chi_k) +  \frac{1}{2} \int_{\mathbb{R}^3} \frac{|\chi_k(x)|^2}{|x|} \ dx $$
$$ = \lim_{k \rightarrow \infty} \frac{1}{2} \| \nabla \chi_k\|^2_{L^2}+ \frac{\omega}{2} \| \chi_k\|^2_{L^2} + \frac{1}{16\pi} \| \nabla \varphi_k \|^2_{L^2} = S_\omega^{min} +
 \frac{1}{2} \int_{\mathbb{R}^3} \frac{|\chi_*(x)|^2}{|x|} \ dx.$$

It is well - known that for any sequence $f_k$ in a Hilbert space
$H$ tending weakly (in $H$) to $f_* \in H$, one has
\begin{equation}\label{eq.P1H}
   \liminf_{k \rightarrow \infty} \| f_k \|_H \geq \| f_* \|_H
\end{equation}
and
\begin{equation}\label{eq.P2H}
    \lim_{k \rightarrow \infty} \| f_k - f_*\|_H = 0 \ \ \Longleftrightarrow \ \  \lim_{k \rightarrow \infty} \| f_k\|_H = \| f_*\|_H.
\end{equation}
From \eqref{eq.P1H} we have
$$ S_\omega^{min} +
 \frac{1}{2} \int_{\mathbb{R}^3} \frac{|\chi_*(x)|^2}{|x|} \ dx \geq \| \nabla \chi_*\|^2_{L^2}+ \frac{\omega}{2} \| \chi_*\|^2_{L^2} + \frac{1}{16\pi} \| \nabla \varphi_* \|^2_{L^2} $$
 and a strict inequality is impossible since this will contradicts the definition of $ S_\omega^{min}.$ Hence
 $$ \lim_{k \rightarrow \infty} \frac{1}{2} \| \nabla \chi_k\|^2_{L^2}+ \frac{\omega}{2} \| \chi_k\|^2_{L^2} + \frac{1}{16\pi} \| \nabla \varphi_k \|^2_{L^2} = $$
 $$ =\frac{1}{2}  \| \nabla \chi_*\|^2_{L^2}+ \frac{\omega}{2} \| \chi_*\|^2_{L^2} + \frac{1}{16\pi} \| \nabla \varphi_* \|^2_{L^2} $$
 and applying \eqref{eq.P2H} we conclude that
 $$ \lim_{k \rightarrow \infty} \|\chi_k - \chi_*\|_{H^1} = 0. $$
 This completes the proof of the Theorem.

 \end{pf3}

\section{Connection between the action and energy minimization
problems}\label{app.C}

Consider the minimization problem
\begin{equation}\label{act-minaa}
    S_{\omega}^{min} = \min \{ S_{\omega}(\chi) ; \chi \in H^1\},
\end{equation}
associated with the action functional $S_{\omega}(\chi)$ and the
Lions--Cazenave minimization problem
\begin{equation}\label{enrg-min}
    I_N = \min\{ \mathcal{E}(\chi) ; \chi \in H^1, \|\chi\|^2_{L^2} = N
    \}.
\end{equation}

As we have seen before, for every $\omega\in (1/16,1/4)$, there
exists (at most one) solution $\chi_{\omega}\in
H^1(\mathbb{R}^3)$, which is positive and radially symmetric, and
such that
\begin{equation}\label{Somega}
    S_{\omega}(\chi_{\omega})=S_{\omega}^{min} .
\end{equation}

Let us denote
\begin{equation}\label{Nomega}
    N(\omega)=\|\chi_{\omega}\|^2_{L^2}.
\end{equation}

The above definition of the function $N(\omega)$ poses the question if
$$ S_{\omega}^{min} = I_{N(\omega)}+\frac{\omega}{2}N(\omega).$$
 For completeness, in this section we shall
 prove the following Lemma.

\begin{lem}\label{Equiv} If $\chi_1$ is a solution of \eqref{enrg-min} with $ N=N(\omega)$, then
$\chi_1$ satisfies the equation
\begin{eqnarray}\label{eq.Hart1}
 -\Delta \chi_1(x) + \int_{\mathbb{R}^3}\frac{\chi_1^2(y) dy}{|x-y|} \ \chi_1(x) - \frac{\chi_1(x)}{|x|} + \omega \chi_1(x) =  0.
\end{eqnarray}
and $$S_\omega(\chi_1) = \min \{ S_{\omega}(\chi) ; \chi \in H^1\}.$$
\end{lem}
\begin{pf} To prove the Lemma we shall follow the idea of the proof of  Corollary 8.3.8 in \cite{Caz}.
It is obvious, that the relation
$$S_{\omega}(\chi_1)=\mathcal{E}(\chi_1 )+\frac{\omega}{2}N(\omega),$$
guarantees that $\chi_1$ is a minimizer of the problem
$$\min_{\|\chi\|^2_{L^2}=N(\omega)} S_{\omega}(\chi).$$
Since,
$$S_\omega(\chi_1) = \min_{\|\chi\|^2_{L^2}=N(\omega)} S_{\omega}(\chi) \geq \min S_{\omega}(\chi) = S_\omega(\chi_\omega), $$
we can use \eqref{Nomega} and see that this inequality becomes equality,
so
 $$S_\omega(\chi_1) = \min_{\|\chi\|^2_{L^2}=N(\omega)} S_{\omega}(\chi) = \min S_{\omega}(\chi) = S_\omega(\chi_\omega). $$
Now, the uniqueness result of Theorem  \ref{th:unique} implies
$\chi_1=\chi_\omega$ and completes the proof.

\end{pf}

\bibliographystyle{model1a-num-names}

\end{document}